\documentclass[a4paper,10pt]{article}
\usepackage[latin1]{inputenc}
\usepackage[T1]{fontenc}
\usepackage{amsmath}
\usepackage{amsfonts}
\usepackage{amstext}
\usepackage{amsthm}
\usepackage[all]{xy}
\usepackage{geometry}
\usepackage{srcltx}
\usepackage{graphics}
\usepackage{epsfig}
\usepackage{subfigure}
\usepackage{slashed}
\usepackage{color}
\usepackage{amssymb}
\usepackage{srcltx}
\newtheorem{theorem}{Theorem}[section]

\newtheorem{rem}[theorem]{Remark}
\newtheorem{prop}[theorem]{Proposition}

\DeclareMathOperator{\im}{im}

\DeclareMathOperator{\sgn}{sgn}

\title{Yet another proof of the Morse index theorem}
\author{Alessandro Portaluri and Nils Waterstraat}

\begin{document}
\date{}
\maketitle

\footnotetext[2]{A. Portaluri was supported by the grant PRIN2009 ``Critical Point Theory and Perturbative Methods
for Nonlinear Differential Equations''.}
\footnotetext[3]{N. Waterstraat was supported by the Berlin Mathematical School and the SFB 647 ``Space--Time--Matter''.}

\begin{abstract}
We give a new analytical proof of the Morse index theorem for geodesics in Riemannian manifolds. 
\end{abstract}

\section{Introduction}
Let $(M,g)$ be a Riemannian manifold of dimension $n\in\mathbb{N}$, Levi-Civita connection $\nabla$ and curvature $R$. Following Klingenberg's monograph \cite{Klingenberg}, we denote by $\Omega_{pq}$ the space of all $H^1$-paths in $M$ joining two given points $p\neq q\in M$, and we recall that $\Omega_{pq}$ carries a natural differentiable structure of a Hilbert manifold such that the energy functional

\begin{align}\label{energy}
E:\Omega_{pq}\rightarrow\mathbb{R},\quad E(\gamma)=\frac{1}{2}\int^1_0{g(\gamma'(x),\gamma'(x))\, dx}
\end{align}
is smooth. Since the geodesics $\gamma:I\rightarrow M$ joining $p$ and $q$ are precisely the critical points of $E$, they play an important role for studying the topology of $\Omega_{pq}$ by Morse-theoretic methods. Here and subsequently, $I$ denotes the unit interval $[0,1]$. Given a critical point $\gamma\in\Omega_{pq}$, we can build the associated Hessian $D^2_\gamma E$ of $E$ at $\gamma$, which is a bilinear form on the tangent space $T_\gamma\Omega_{pq}$ at $\gamma$. The index $\mu(\gamma)$ is by definition the Morse index of $D^2_\gamma E$, that is, the dimension of a maximal subspace of $T_\gamma\Omega_{pq}$ on which the quadratic form induced by $D^2_\gamma E$ is negative definite. The fundamental theorem of Morse theory states that if $M$ is complete and $D^2_\gamma E$ is non-degenerate at all critical points, then $\Omega_{pq}$ has the homotopy type of a countable CW-complex which contains one cell of dimension $\lambda$ for each geodesic from $p$ to $q$ of index $\lambda$ (cf. \cite[Theorem 17.3]{MilnorMorse}). Consequently, it is important to know the Morse index of a geodesic, and the famous Morse index theorem asserts that $\mu(\gamma)$ is just the total number of conjugate points along $\gamma$. This allows to compute the Morse index for geodesics in various types of manifolds $M$ (cf. eg. Part IV of Milnor's monograph \cite{MilnorMorse}).\\ 
Nowadays several proofs of the Morse index theorem are known. The oldest reference of which we are aware is Morse's classical work \cite{MorseTrans} (cf. also \cite{Morse}). Different proofs were provided later on, for example by Ambrose \cite{Ambrose}, Duistermaat \cite{Duistermaat} and Osborn \cite{OsbornI}, \cite{OsbornII}. Smale proved in \cite{Smale}, \cite{SmaleCorr} a Morse index theorem for strongly elliptic differential operators acting on Euclidean bundles over Riemannian manifolds with boundary, which reduces to the assertion for geodesics if the underlying manifold is the unit interval. Later Uhlenbeck \cite{Uhlenbeck} and Swanson \cite{Swansona},\cite{Swansonb} developed new proofs of Smale's theorem. Most of these proofs of the Morse index theorem either use domain monotonicity properties of eigenvalues for strongly elliptic operators or intersection theory in Grassmannian manifolds.\\
Another direction of research around the Morse index theorem was initiated by Helfer in \cite{HELFER}, who studied geodesics in semi-Riemannian manifolds and observed that the number of conjugate points along geodesics can now be infinite. Consequently, the classical Morse index theorem does not extend to semi-Riemannian manifolds. However, Helfer proved a generalisation which in particular covers the classical result for Riemannian manifolds. General Morse index theorems for geodesics in semi-Riemannian manifolds were later proved by Piccione and Tausk (cf. \cite{PiccioneMITinSRG} and the references therein), and by the first author in collaboration with Musso and Pejsachowicz \cite{MPP} (cf. also \cite{Wa} for a $K$-theoretic proof of it). Of course, also these theorems reduce to the classical one for geodesics in Riemannian manifolds. In the latter case this follows from \cite[Prop. 3.4]{MPP} and the slightly hidden computation in the last paragraph of page 97 in \cite{MPP}, if one uses the fact that conjugate points are now isolated (cf. Prop. 2.2 in \cite{Mercuri}).\\
The aim of this article is to introduce a new proof of the the Riemannian Morse index theorem, which combines certain ideas from \cite{MPP} with arguments recently obtained by the authors in \cite{AleIchDomain} and \cite{AleIchBalls}, respectively. The outcome is a proof which is based on an essentially well known assertion about the Morse index of quadratic forms in Hilbert spaces coming from variational bifurcation theory \cite{FPR} (cf. also \cite[Thm. 7.3]{Rabier}). We recall this fact in Section \ref{section-Morse} in order to make our presentation self-contained, particularly for differential geometers. Subsequently, the proof of the equality of $\mu(\gamma)$ and the total number of conjugate points along $\gamma$ boils down to an elementary computation which we show in Section \ref{section-proof}. Let us point out that our proof neither uses the domain monotonicity of eigenvalues of strongly elliptic equations nor arguments from intersection theory.


\section{The Morse index theorem}
As before, we follow the presentation of Klingenberg's monograph \cite{Klingenberg}. The Sobolev space $H^1(I,\mathbb{R}^n)$ of all absolutely continuous curves having a square integrable derivative is a Hilbert space with respect to the scalar product

\begin{align}\label{scalarprod}
\langle \gamma_1,\gamma_2\rangle_{H^1(I,\mathbb{R}^n)}=\langle\gamma_1,\gamma_2\rangle_{L^2(I,\mathbb{R}^n)}+\langle\gamma'_1,\gamma'_2\rangle_{L^2(I,\mathbb{R}^n)}.
\end{align}
Let $H^1(I,M)$ be the set of all curves $\gamma:I\rightarrow M$ which belong to $H^1(I,\mathbb{R}^n)$ with respect to charts. There is a natural smooth atlas on $H^1(I,M)$ having $H^1(I,\mathbb{R}^n)$ as model space, where the charts are defined by means of the exponential map on $M$. The tangent space $T_\gamma H^1(I,M)$ at an element $\gamma\in H^1(I,M)$ can be identified canonically with the space of all $H^1$-curves along $\gamma$, that is, the space of all $\eta\in H^1(I,TM)$ such that  $\pi\circ\eta=\gamma$, where $\pi:TM\rightarrow M$ denotes the projection of the tangent bundle of $M$.\\
Let us now consider two points $p\neq q\in M$. Since the endpoint evaluation

\[P:H^1(I,M)\rightarrow M\times M,\quad \gamma\mapsto(\gamma(0),\gamma(1))\]
is a submersion, it follows that $\Omega_{pq}:=P^{-1}(\{p,q\})$, which is the set of all curves joining $p$ and $q$, is a $2n$-codimensional smooth submanifold of $H^1(I,M)$ having the model space

\[H^1_0(I,\mathbb{R}^{n})=\{\gamma\in H^1(I,\mathbb{R}^n):\, \gamma(0)=\gamma(1)=0\}.\] 
Moreover, the tangent space $T_\gamma\Omega_{pq}$ can be identified with the space of all $H^1$ vector fields along $\gamma$ that vanish at $\gamma(0)$ and $\gamma(1)$. The energy functional \eqref{energy} is a smooth map on $\Omega_{pq}$ and the critical points of $E$ are the geodesics between $p$ and $q$.  Let now $\gamma\in\Omega_{pq}$ be a geodesic. The Hessian at the critical point $\gamma$ of $E$ is the bilinear form 

\[D^2_\gamma E:T_\gamma \Omega_{pq}\times T_\gamma \Omega_{pq}\rightarrow\mathbb{R}\]
given by  

\begin{align}\label{Hessian}
D^2_\gamma E(\xi,\eta)=\int^1_0{g\left(\frac{\nabla}{dx}\xi(x),\frac{\nabla}{dx}\eta(x)\right)\,dx}-\int^1_0{g(R(\gamma'(x),\xi(x))\gamma'(x),\eta(x))\, dx}.
\end{align}
We denote by $\mu(\gamma)$ the Morse index of $\gamma$, which is the dimension of the maximal subspace of $T_\gamma \Omega_{pq}$ on which the quadratic form induced by $D^2_\gamma E$ is negative definite.\\
Let now $\Gamma(\gamma)$ be the space of all smooth vector fields along $\gamma$, and consider the Jacobi equation

\begin{align}\label{Jacobi}
\frac{\nabla^2}{dx^2}\xi(x)+R(\gamma'(x),\xi(x))\gamma'(x)=0,\quad x\in I,
\end{align}
for vector fields $\xi\in\Gamma(\gamma)$. Let us recall that an instant $t\in I$ is said to be conjugate if 

\[m(t)=\dim\{\xi\in\Gamma(\gamma):\, \xi\,\,\text{satisfies}\, \eqref{Jacobi},\, \xi(0)=0,\,\xi(t)=0\}>0. \]
Finally, $\gamma$ is called non-degenerate if $m(1)=0$. With all this said, we now can state the well known Morse Index Theorem:

\begin{theorem}\label{MIT}
Let $\gamma$ be a non-degenerate geodesic. Then $m(t)=0$ for all but a finite number of $t\in I$, and

\[\mu(\gamma)=\sum_{t\in I}{m(t)}.\]
\end{theorem}


\section{The proof}

\subsection{Morse index and crossing forms}\label{section-Morse}
The aim of this section is to prove an abstract result about the Morse index of quadratic forms on Hilbert spaces, which will turn out to be more than half the battle for obtaining the Morse Index Theorem \ref{MIT}.\\
Concretely, let $H$ be a real Hilbert space and $q:H\rightarrow\mathbb{R}$ a bounded quadratic form, that is, there exists a bounded bilinear form $b:H\times H\rightarrow\mathbb{R}$ such that $q(u)=b(u,u)$, $u\in H$. We will assume throughout that the Riesz representation of $b$, which is the unique selfadjoint operator $L:H\rightarrow H$ that satisfies

\[b(u,v)=\langle Lu,v\rangle,\quad u,v\in H,\]
is of the form $L=I_H+K$, where $K:H\rightarrow H$ is compact and $I_H$ denotes the identity on $H$. Under this assumption, the Morse index of $q$,

\begin{align}\label{finite}
\mu(q)=\sup\dim\{U\subset H:q(u)<0,\, u\in U\setminus\{0\}\},
\end{align}
is finite and coincides with the number of negative eigenvalues of $L$ counted according to their multiplicities. In what follows we will also denote the number \eqref{finite} by $\mu(L)$. For later reference, let us point out the following three elementary properties of the Morse index:

\begin{enumerate}
	\item[i)] If $U:H\rightarrow H$ is an orthogonal operator, then $\mu(U^{-1} LU)=\mu(L)$.
	\item[ii)] If $L$ is reduced by a splitting $H=H_0\oplus H_1$, that is, $L(H_i)\subset H_i$, $i=0,1$, then 
	
	\[\mu(L)=\mu(L\mid_{H_0})+\mu(L\mid_{H_1}).\]
	\item[iii)] If $L_0,L_1$ belong to the same component of 
	
\[GL^\textup{sa}_c(H)=\{I_H+K\in GL(H):\, K\,\text{compact, selfadjoint}\},\]
then $\mu(L_1)=\mu(L_0)$. 
\end{enumerate}

Let us now assume that we can connect $K_1$ to an operator $K_0$ through a $C^2$ path $K_\lambda$, $\lambda\in I$, of compact selfadjoint operators. This induces a path of quadratic forms by

\[q:I\times H\rightarrow\mathbb{R},\quad q_\lambda(u)=\langle L_\lambda u,u\rangle,\]
where $L_\lambda=I_H+K_\lambda$. We call an instant $\lambda\in I$ a crossing of $q$ if $\ker L_\lambda\neq\{0\}$,
or equivalently, if $q_\lambda$ is degenerate. Moreover, we call a crossing $\lambda$ regular if the restriction of the quadratic form $\dot q_\lambda:H\rightarrow\mathbb{R}$ to $\ker L_\lambda$ is non-degenerate, where $\dot{}$ denotes the derivative with respect to the parameter $\lambda$. Note that $\dot q_\lambda(u)=\langle\dot K_\lambda u,u\rangle$, $u\in H$. The abstract result on the Morse index now reads as follows:

\begin{prop}\label{prop-crossing}
If $q_0,q_1$ are non-degenerate and all crossings of $q$ in $I$ are regular, then there are only finitely many crossings of $q$ and

\begin{align}\label{crossing}
\mu(q_0)-\mu(q_1)=\sum_{\lambda\in I}{\sgn\left(\dot{q}_\lambda\mid_{\ker L_\lambda}\right)}.
\end{align}
\end{prop}

In our proof of the Morse Index Theorem below we will join the quadratic form induced by \eqref{Hessian} to a positive definite one, and consequently the right hand side of \eqref{crossing} gives the Morse index.

\begin{proof}
That a regular crossing $\lambda_0$ of $q$ is isolated clearly follows from the assertion that there exist positive numbers $\varepsilon, c>0$ such that

\begin{align}\label{estimate}
\|L_\lambda u\|\geq c\,|\lambda-\lambda_0|\|u\|,\quad |\lambda-\lambda_0|<\varepsilon,\,\, u\in H,
\end{align}
which we are going to prove now. Since $L_{\lambda_0}$ is selfadjoint and has a closed range, we have an orthogonal decomposition $H=\ker L_{\lambda_0}\oplus\im L_{\lambda_0}$. In what follows, we set for notational convenience $H_0=\ker L_{\lambda_0}$. That \eqref{estimate} holds for all $u\in H_0$ and all $\lambda$ which are sufficiently close to $\lambda_0$ is an immediate consequence of the assumption that $\dot q_{\lambda_0}$ is non-degenerate on $H_0$ and we leave the details to the reader. Let now $u\in\im L_{\lambda_0}=H^\perp_0$ and let us denote by $P$ the orthogonal projection onto the closed subspace $H^\perp_0\subset H$. Since $GL(H^\perp_0)$ is open in the Banach space $\mathcal{L}(H^\perp_0)$ of all bounded operators on $H^\perp_0$, $PL_\lambda P:H^\perp_0\rightarrow H^\perp_0$ is an isomorphism for all $\lambda$ sufficiently close to $\lambda_0$. Consequently, there exists $\varepsilon>0$ and $c>0$ such that

\[\|PL_\lambda Pu\|\geq c\,\|u\|,\quad u\in H^\perp_0,\, |\lambda-\lambda_0|<\varepsilon.\]
From $Pu=u$ for $u\in H^\perp_0$ and $\|PL_\lambda u\|\leq\|P\|\|L_\lambda u\|=\|L_\lambda u\|$, we obtain $\|L_\lambda u\|\geq c\|u\|$ for all $|\lambda-\lambda_0|<\varepsilon$ and $u\in H^\perp_0$ and this gives \eqref{estimate} when we assume that $\varepsilon<1$.\\
Consequently, regular crossings are isolated, and by iii) we can henceforth assume that there is only a single regular crossing $\lambda_0$ of $q$. Then $\ker L_{\lambda_0}\neq\{0\}$, and since $L_{\lambda_0}$ is a compact perturbation of the identity, there exists $\varepsilon>0$ such that $0$ is the only eigenvalue of $L_{\lambda_0}$ in the interval $[-\varepsilon,\varepsilon]$. Moreover, by the continuity of finite parts of spectra (cf. \cite[Theorem I.II. 4.2]{GohbergClasses}), there exists $\rho>0$ such that $\pm\varepsilon$ is not in the spectrum of $L_\lambda$ for all $\lambda\in[\lambda_0-\rho,\lambda_0+\rho]$. Let now $P_\lambda$ denote the orthogonal projection onto the sum of the eigenspaces of $L_\lambda$ with respect to eigenvalues in $[-\varepsilon,\varepsilon]$. Then $P$ is a continuously differentiable path of bounded projections, and according to \cite[Sect. VI.2]{Kato}, there exists an interval $[a,b]$, $a<b$, containing $\lambda_0$ and a continuously differentiable path $U_\lambda$ of orthogonal operators such that

\begin{align}\label{orthogonal}
U_{\lambda_0}=I_H,\quad U_\lambda P_\lambda U^{-1}_\lambda=P_{\lambda_0},\quad\lambda\in[a,b].
\end{align}  
Since $L_\lambda$ is reduced by the decomposition $H=\ker P_\lambda\oplus\im P_\lambda$, we deduce from \eqref{orthogonal} that the operator $U_\lambda L_\lambda U^{-1}_\lambda$ is reduced by $H=H_0\oplus H^\perp_0$. Moreover, the restriction $U_\lambda L_\lambda U^{-1}_\lambda\mid_{H^\perp_0}:H^\perp_0\rightarrow H^\perp_0$ is an isomorphism. We conclude from the properties i)-iii) of the Morse index

\begin{align}\label{step1}
\begin{split}
\mu(q_a)-\mu(q_b)&=\mu(L_a)-\mu(L_b)=\mu(U_aL_aU^{-1}_a)-\mu(U_bL_bU^{-1}_b)\\
&=\mu(U_aL_aU^{-1}_a\mid_{H_0})+\mu(U_aL_aU^{-1}_a\mid_{H^\perp_0})\\
&-(\mu(U_bL_bU^{-1}_b\mid_{H_0})+\mu(U_bL_bU^{-1}_b\mid_{H^\perp_0}))\\
&=\mu(U_aL_aU^{-1}_a\mid_{H_0})-\mu(U_bL_bU^{-1}_b\mid_{H_0}).
\end{split}
\end{align}
Consequently, let us now consider the path of symmetric operators $\ell_\lambda=U_\lambda L_\lambda U^{-1}_\lambda\mid_{H_0}$, $\lambda\in[a,b]$, on the finite dimensional space $H_0$. We obtain for $\lambda\in[a,b]$

\begin{align*}
\dot\ell_\lambda=\frac{d}{d\lambda}(U_\lambda)\circ L_\lambda\circ U^{-1}_\lambda\mid_{H_0}+U_\lambda\circ\frac{d}{d\lambda}(L_\lambda)\circ U^{-1}_\lambda\mid_{H_0}+U_\lambda\circ L_\lambda\circ\frac{d}{d\lambda}(U^{-1}_\lambda)\mid_{H_0}
\end{align*}
and, using that $U_{\lambda_0}=I_H$, we conclude for $u,v\in H_0$

\begin{align*}
\langle\dot\ell_{\lambda_0}u,v\rangle&=\langle\dot L_{\lambda_0}u,v\rangle+\langle  L_{\lambda_0}\dot U^{-1}_{\lambda_0}u,v\rangle=\langle\dot L_{\lambda_0}u,v\rangle+\langle\dot U^{-1}_{\lambda_0}u,L_{\lambda_0}v\rangle\\
&=\langle\dot L_{\lambda_0}u,v\rangle.
\end{align*}
Hence, 

\[\dot q_{\lambda_0}(u)=\langle\dot\ell_{\lambda_0}u,u\rangle,\quad u\in\ker L_{\lambda_0},\]
and we note that by \eqref{step1} the equality \eqref{crossing} is shown if we prove that 

\begin{align}\label{step2}
\mu(\ell_a)-\mu(\ell_b)=\sgn\langle\dot\ell_{\lambda_0}\cdot,\cdot\rangle.
\end{align}
Since $\dot q_{\lambda_0}:H_0\rightarrow\mathbb{R}$ is non-degenerate and $\dot q_{\lambda_0}(u)=\langle\dot\ell_{\lambda_0}u,u\rangle$, $u\in H_0$, we infer that $\dot\ell_{\lambda_0}$ is an invertible operator on the finite dimensional space $H_0$. Hence, there exists $\alpha>0$ such that $\dot\ell_{\lambda_0}+A$ is invertible for all linear operators $A:H_0\rightarrow H_0$ of norm less than $\alpha$. We now consider for $t\in[0,1]$ and $\lambda=a,b$

\[H(t,\lambda)=t\ell_\lambda+(1-t)(\lambda-\lambda_0)\dot\ell_{\lambda_0}=(\lambda-\lambda_0)\left(t\left(\frac{\ell_\lambda-\ell_{\lambda_0}}{\lambda-\lambda_0}-\dot\ell_{\lambda_0}\right)+\dot\ell_{\lambda_0}\right),\]
where we use that by definition $\ell_{\lambda_0}=L_{\lambda_0}\mid_{H_0}=0$. Since we may choose $a<\lambda_0$ and $b>\lambda_0$ in \eqref{orthogonal} arbitrarily close to $\lambda_0$, we can assume that the norm of $\frac{\ell_\lambda-\ell_{\lambda_0}}{\lambda-\lambda_0}-\dot\ell_{\lambda_0}$ is less than $\alpha$ for $\lambda=a,b$. Consequently, $H(\lambda,t)$ is invertible for $t\in[0,1]$ and $\lambda=a,b$, and we conclude by iii) that

\[\mu(\ell_a)=\mu((a-\lambda_0)\,\dot\ell_{\lambda_0})=\mu(-\dot\ell_{\lambda_0})\quad\text{and}\quad \mu(\ell_b)=\mu((b-\lambda_0)\,\dot\ell_{\lambda_0})=\mu(\dot\ell_{\lambda_0}).\]
This shows \eqref{step2} since 
\[\sgn\langle\dot\ell_{\lambda_0}\cdot,\cdot\rangle=\mu(-\dot\ell_{\lambda_0})-\mu(\dot\ell_{\lambda_0})\]
by the very definition of the signature.
\end{proof}

\begin{rem}
A more general version of Proposition \ref{prop-crossing} can be found in \cite[Thm 4.1]{FPR}, where the operators $L_\lambda$, $\lambda\in I$, are selfadjoint Fredholm operators. Then the difference of the Morse indices in \eqref{crossing} is no longer defined in general. Instead, the left hand side in \eqref{crossing} is given by the spectral flow, which is a homotopy invariant for paths of selfadjoint Fredholm operators that, in contrast to the situation in Proposition \ref{prop-crossing}, generally depends on the whole path $q$ and not only on its endpoints. 
\end{rem}


\subsection{Proof of the Morse index theorem}\label{section-proof}
We now consider the Hessian \eqref{Hessian} and the Jacobi equation \eqref{Jacobi} of a given geodesic $\gamma\in\Omega_{pq}$, where $p\neq q$ are two points in $M$. We take an orthonormal frame along $\gamma$ made by $n$ parallel vector fields $\{e^1,\ldots,e^n\}$. This allows us to write any vector field $\xi$ along $\gamma$ uniquely as

\begin{align}\label{baserep}
\xi(x)=\sum^n_{i=1}{u_i(x)\,e^i(x)}, \quad x\in I,
\end{align}
and yields an identification of vector fields along $\gamma$ that belong to $T_\gamma\Omega_{pq}$ and $H^1_0(I,\mathbb{R}^n)$. In what follows, we use on $H^1_0(I,\mathbb{R}^n)$ the scalar product

\begin{align}\label{scalarprodII} 
\langle u,v\rangle_{H^1_0(I,\mathbb{R}^n)}=\langle u',v'\rangle_{L^2(I,\mathbb{R}^n)},\quad u,v\in H^1_0(I,\mathbb{R}^n),
\end{align}
which is easily seen to be equivalent to the one induced by \eqref{scalarprod}.\\
Plugging \eqref{baserep} into the quadratic form of the Hessian \eqref{Hessian} gives

\begin{align*}
q_1(u)=\int^1_0{\langle u'(x),u'(x)\rangle\,dx}-\int^1_0{\langle S(x)u(x),u(x)\rangle\,dx},
\end{align*}
where $S$ denotes the smooth path of symmetric matrices having components 

\[S_{ij}(x)=g(R(\gamma'(x),e^i(x))\gamma'(x),e^j(x)),\quad x\in I,\,\,1\leq i,j\leq n.\]
Moreover, when plugging \eqref{baserep} into \eqref{Jacobi}, the Jacobi equation transforms to

\begin{align}\label{JacobiII}
u''(x)+S(x)u(x)=0,\quad x\in I.
\end{align}
We now join $q_1$ to a positive definite quadratic form on $H^1_0(I,\mathbb{R}^n)$ by the path

\begin{align*}
q_\lambda(u)=\int^1_0{\langle u'(x),u'(x)\rangle\,dx}-\int^1_0{\langle \lambda^2S(\lambda\cdot x)u(x),u(x)\rangle\,dx},\quad \lambda\in I.
\end{align*} 
Note that $q_0$ is just the quadratic form induced by the scalar product \eqref{scalarprodII} on $H^1_0(I,\mathbb{R}^n)$. From the compactness of the inclusion $H^1_0(I,\mathbb{R}^n)\hookrightarrow C(I,\mathbb{R}^n)$, it is readily seen that the Riesz representation of $q_\lambda$ is of the form $I_{H^1_0(I,\mathbb{R}^n)}+K_\lambda$ for some smooth path $K_\lambda$ of compact operators on $H^1_0(I,\mathbb{R}^n)$. Finally, since $K_0=0$ and $q_1$ is non-degenerate by assumption, neither $0$ nor $1$ is a crossing of the path $q$.\\
In view of Proposition \ref{prop-crossing}, our aim is to show that each crossing $\lambda_0$ of $q_\lambda$ is regular and that $\sgn\dot q_{\lambda_0}\mid_{\ker L_{\lambda_0}}=-m(\lambda_0)$. Let $\lambda_0\in(0,1)$ be a crossing of $q$ and let $u\in\ker L_{\lambda_0}\subset H^1_0(I,\mathbb{R}^n)$. It is clear that $u$ is smooth by elliptic regularity, and moreover, $u$ solves the boundary value problem

\begin{equation}\label{bvp}
\left\{
\begin{aligned}
\frac{d^2u}{dx^2}+\lambda^2_0\,S(\lambda_0\cdot x)u(x)&=0,\,\,x\in I,\\
u(0)=u(1)&=0.
\end{aligned}
\right.
\end{equation}
We introduce a family of smooth functions by $u^\lambda_{\lambda_0}(x)=u(\frac{\lambda}{\lambda_0}x)$, $0\leq\frac{\lambda}{\lambda_0}x\leq 1$, and we denote henceforth 

\[\dot u(x)=\frac{d}{d\lambda}\mid_{\lambda=\lambda_0}u^\lambda_{\lambda_0}(x)=\frac{x}{\lambda_0}\frac{du}{dx},\quad x\in(0,1).\]
From \eqref{bvp} it is readily seen that

\begin{align}\label{equation}
\frac{d^2}{dx^2}u^\lambda_{\lambda_0}(x)+\lambda^2S(\lambda\cdot x)u^\lambda_{\lambda_0}(x)=0,\quad 0<\frac{\lambda}{\lambda_0}x< 1, 
\end{align}
and by differentiating \eqref{equation} with respect to $\lambda$ and evaluating at $\lambda=\lambda_0$, we obtain

\[\frac{d^2\dot u}{dx^2}+\frac{d}{d\lambda}\mid_{\lambda=\lambda_0}(\lambda^2S(\lambda\cdot x))u(x)+\lambda^2_0S(\lambda_0\,x)\dot u(x)=0,\quad x\in(0,1).\]
We take scalar products with $u$, integrate and obtain

\[\int^1_0{\langle\frac{d^2\dot u}{dx^2},u(x)\rangle\, dx}+\int^1_0{\langle \frac{d}{d\lambda}\mid_{\lambda=\lambda_0}(\lambda^2S(\lambda\cdot x))u(x),u(x)\rangle\,dx}+\int^1_0{\langle \lambda^2_0S(\lambda_0\,x)\dot u(x),u(x)\rangle\, dx}=0.\]
Performing integration by parts twice and using \eqref{bvp}, we finally conclude that

\[\dot q_{\lambda_0}(u)=-\frac{d}{d\lambda}\mid_{\lambda=\lambda_0}\int^1_0{\lambda^2 \langle S(\lambda\cdot x)u(x),u(x)\rangle\,dx}=-\langle \dot u(1),\frac{du}{dx}(1)\rangle+\langle\dot u(0),\frac{d u}{dx}(0)\rangle =-\frac{1}{\lambda_0}\|\frac{du}{dx}(1)\|^2.\]
Since $u$ solves \eqref{bvp} and $u\neq 0$, $\frac{du}{dx}(1)$ cannot vanish.\\
Consequently, we have shown that $\dot q_{\lambda_0}$ is negative definite on $\ker L_{\lambda_0}$. This means in particular that $\dot q_{\lambda_0}\mid_{\ker L_{\lambda_0}}$ is non-degenerate and $\sgn\dot q_{\lambda_0}\mid_{\ker L_{\lambda_0}}=-\dim\ker L_{\lambda_0}$. Finally, elements in the kernel of $q_{\lambda_0}$ are exactly the solutions of \eqref{bvp}, which are in turn just the solutions of \eqref{JacobiII} that satisfy $u(0)=u(\lambda_0)=0$. This shows that $\dim\ker L_{\lambda_0}=m(\lambda_0)$, and now the Morse Index Theorem \ref{MIT} follows from Proposition \ref{prop-crossing}.


\vspace{1cm}
Alessandro Portaluri\\
Department of Agriculture, Forest and Food Sciences\\
Universit\`a degli studi di Torino\\
Largo Paolo Braccini, 2\\
10095 Grugliasco (TO)\\
Italy\\
E-mail: alessandro.portaluri@unito.it

\vspace{1cm}
Nils Waterstraat\\
Institut für Mathematik\\
Humboldt-Universität zu Berlin\\
Unter den Linden 6\\
10099 Berlin\\
Germany\\
E-mail: waterstn@math.hu-berlin.de

\end{document}